\documentclass[a4paper,twoside,reqno]{amsart}

\usepackage{url}
\usepackage{amssymb,amscd}
\usepackage[mathscr]{eucal}

\newcommand{\FF}{\mathbb{F}}

\newcommand{\LL}{\mathcal{L}}

\theoremstyle{plain}
\newtheorem{prop}{Proposition}

\theoremstyle{definition}

\newtheorem{rem}[prop]{Remark}

\theoremstyle{remark}

\newtheorem{Examples}[prop]{Examples}
\newtheorem{Rems}[prop]{Remarks}

\newcounter{reml}

\newenvironment{remlist}{\begingroup\setcounter{reml}{0}}
  {\endgroup}

\title[Hyperovals in $Q^+(6,4)$]
{On hyperovals in $Q^+(6,4)$}
\author{Dmitrii Pasechnik}
\address{Department of Computer Science, Oxford University, UK}

\begin{document}
\begin{abstract}
According to a computer search described in \cite{Pa:epolsp},
in $Q^+(6,4)$ there are two types of hyperovals, having 72 and 96 points,
respectively. Here we give geometric descriptions for these examples.
\end{abstract}

\maketitle

\section{Introduction}
A hyperoval in a partial linear space is a subset of points intersecting each line in either 0 or 2 points.
Classical examples are hyperovals in projective planes of order $2^k$. In particular, hyperovals
in the projective plane $PG(2,4)$ over $\FF_4$ (i.e. of order 4)
appear as building blocks of Witt designs, leading to Mathieu sporadic simple
groups $M_{22}$, $M_{23}$ and $M_{24}$. More generally,
hyperovals in polar spaces over $\FF_4$ lead to more sporadic simple groups, these of Fischer, see
\cite{Pa:epolsp}, where this has been investigated, in part relying on computer searchers, and
further, to Baby Monster, see \cite{MR1827737}.

One of these searches in \cite{Pa:epolsp}
succeeded in enumerating the hyperovals in $Q^+(6,4)$, the line Grassmanian $\LL$ of $PG(3,4)$. 
There are two examples (up to the group action),
on 72 and on 96 points.

In this note we provide a geometric interpretation of the examples of hyperovals in $\LL$ found there.
Other hyperovals studied in \cite{Pa:epolsp}
were further investigated in \cite{MR2658930, MR2935473}.
Combinatorially, hyperovals in $\LL$ are locally $5\times 5$-grid graphs, recently studied in \cite{2019arXiv190307931A},
a particular type of extended
generalised quadrangles, see e.g. \cite{CHP,Pas:book}.

Recall that the \emph{lines} of the Grassmannian $\LL$ consist of the 5 lines of $\Pi$ through a point $p$ on a plane
$P$. We will refer to them as \emph{pencils} and denote them by $(p,P)$, for $p\in P$, to avoid confusion between
the lines of  $\Pi:=PG(3,4)$ an the lines of $\LL$. That is, a hyperoval of $\LL$ is a set of lines
of $\Pi$ that intersects each pencil in 0 or 2 lines.

\section{Geometric constructions}
\paragraph{\bf The 72-point example}
Let $H$ be a hyperbolic quadric in $\Pi$, that is, the quadric of $+$ type, with the
automorphism group $PGO^+_4(4)$. There are two classes, each of size 5, of mutually skew lines
on $H$; together they form a $5\times 5$ grid, each of the classes covers the 25 points of $\Pi$ on $H$, and,
dually, there are 25 planes of $\Pi$ intersecting $H$ in 9 points (which lie on the union of two intersecting lines
on $H$). See e.g. \cite[Sect.~15.3]{MR840877} for details. The remaining 60 planes of $\Pi$ intersect $H$ in the 5 point of a conic; dually, each point of $\Pi$ not
on $H$ lies on 5 planes intersecting $H$ in 9 points.
Out of 357 lines of $\Pi$, 10 lie on $H$, there is  also a non-empty set $\mathcal{O}$ of lines that do not intersect $H$. 

Each plane $P$ intersecting $H$ in a conic $C$ contains 
6 of lines in $\mathcal{O}$. Indeed, there are 10 lines intersecting $C$ in two points, and 5 intersecting $C$ in one point;
the remaining 6 do not intersect $C$, and thus do not interect $H$. We have $60\times 6/5=72=|\mathcal{O}|$.
Dually, each point outside $H$ is on 6 lines from $\mathcal{O}$.

Let $L$ be a pencil $(p,P)$. Let $\ell\in L\cap\mathcal{O}$. Then $P$ intersects $H$ in a conic $C$, and as $p$ is on exactly
two lines in $P$ missing $C$, we see that there is exactly one more line in $L\cap\mathcal{O}$. Hence any pencil of $\LL$
intersects $\mathcal{O}$ in 0 or 2 elements, and we have proved the following.
\begin{prop}\label{ho72}
The $72$ lines skew to a hyperbolic quadric in $\Pi$ form a hyperoval in $\LL$. \qed
\end{prop}

\paragraph{\bf The 96-point example}
Let $\hat{S}$ be a regular line spread (also known as \emph{elliptic congruence}, as it is related to a class
of elliptic quadrics, with $\hat{S}$ tangent lines to any of them) in $\Pi$.
See e.g. \cite[Sect.~17.1]{MR840877} for details.
In particular, $\hat{S}$ consists of
17 lines covering all the 85 points of $\Pi$, and dually, each plane of $\Pi$
contains a line in $\hat{S}$. Let $s\in \hat{S}$, and denote
$S:=\hat{S}\setminus\{s\}$.  The lines of $\Pi$ are partitioned into $S$, the set 
$s^\perp$ - the 101 lines equal to or intersecting $s$, and the set $A$ of
the remaining 240 lines. 
The stabiliser of $\hat{S}$ in $PGL_4(4)$
acts\footnote{This action has a kernel of order 5, acting transitively on the points
of each line.} as $PGL_2(16)=PGO^-_4(4)$ on $\hat{S}$, and the subgroup $G_S$ fixing $s$ in this
action acts as $AGL_2(16)=2^{4}:15$ on $S$. Observe that $G_S$ acts
transitively on $A$; moreover, $G_S$ has a subgroup $G_S^*$ of index 3, which
has 3 orbits $A_1$, $A_2$, $A_3$ on $A$, each of size 80.
With $\zeta\in \mathbb{F}_{16}^*$ of (multiplicative) order 15,
$G_S^*$ may be assumed to be
\[
G_S^*:=\left\langle \begin{pmatrix}\zeta^3 &0\\0&1\end{pmatrix},
\begin{pmatrix}1 &0\\0&\zeta^3\end{pmatrix}, 
\begin{pmatrix}1 &0\\1&1\end{pmatrix}\right\rangle,
\]
and 
$\zeta^3\mapsto \begin{pmatrix}\omega & 1\\\omega^2&1\end{pmatrix}$,
with $\omega\in\mathbb{F}_4^*$ of order 3,  
specifies an embedding of $G_S^*$ into $PSO^{-}_4(4)$.
Note that $|G_S^*|=2^4.5^2=400$.

We are going to show that $\mathcal{O}':=S\cup A_1$ is a hyperoval of $\LL$.

Any pencil $(p,P)$ of $\LL$,
such that either $p\in s$ or $s\in P$, does not contain any element of $S\cup A_1$.
Thus we need to consider the pencils $(p,P)$,
such that $p\not\in s\not\in P$. Such a pencil contains a line $p_s\in s^\perp$ joining $p$ to a point on $s$.

There exists unique $\ell\in S$ in $P$. Let $p\in\ell$. Then the remaining 3 lines of $(p,P)$ lie in one orbit, $A$, of $G_S$.
As 3 does not divide $|G_S^*|$,
these remaining 3 lines lie in different orbits of $G_S^*$. Thus $A_1$ intersects $(p,P)$ in exactly one line $\ell'$, and the 2 lines
in $\mathcal{O}'\cap (p,P)$ are $\ell$ and $\ell'$.

It remains to deal with the case where $p\not\in\ell$.
There are 4 lines in $(p,P)$ which could potentially intersect $\mathcal{O'}$. By the choice of $\ell$, this intersection is contained in $A_1$. 
In the previous case,
we have established that $p_k\in\ell$ is incident to
unique $\ell_k\in A_1$, for $1\leq k\leq 5$.

It turns out that $\ell,\ell_1,\dots,\ell_5$ form a dual hyperoval in
$P$. Indeed, the stabiliser of $P$ in $G_S^*$ is of order 5.
An element of order 5 in the automorphism group of $P$ fixes a point,
which must be $s\cap P$, 
a line not containing this point, which must be $\ell$, and
its orbits on the lines not through the fixed points are three dual
conics, one of which consists of $\ell_1,\dots,\ell_5$.

Therefore $p$ is incident to either exactly 2 lines from
$\ell_1,\dots,\ell_5$, or to none of them. Hence
\begin{prop}\label{ho96}
The $96$ lines in $\mathcal{O}'$ form a hyperoval in $\LL$. \qed
\end{prop}

\begin{rem}
  The full automorphism group of the locally $5\times 5$-grid graph
  associated the 96-line example is 8 times bigger than $G_S^*$.
  In particular, there is an automorphism swapping the two
  classes of the 6-cliques corresponding to points and hyperplanes
  of $\Pi$ not on $s$, and the Galois group of $\mathbb{F}_4$.
\end{rem}

\section{Concluding remarks}

Combining the computer computations and the above observations, one establishes the following.
\begin{prop}\label{prop:nei}
Let $\ell$ be an element of a hyperoval in $\LL$. Then there exists a hyperbolic quadric $H$
in $\Pi$ so that  $\ell^\perp$ in the collinearity graph of $\LL$ coincides
with the $\ell^\perp$ in the $72$-point hyperoval associated with $H$.
\end{prop}
With Proposition~\ref{prop:nei} at hand, it should be is possible to provide a complete classification of
the hyperovals in $\LL$.

The author observed and Antonio Pasini confirmed in a personal
communication that the 72-point construction generalises
for any $PG(3,2^k)$, with $k>2$, providing a diagram geometry
with diagram $o-L-o==o$, where $L$ denotes certain partial
linear space of lines in $PG(2,2^k)$ missing a hyperoval. 

\subsection*{Acknowledgement} We thank Antonio Pasini for a helpful discussion.
The author was supported by EU OpenDreamKit Horizon 2020 project.
GAP \cite{GAP4} was used to carry out various experiments and tests.

\bibliography{geom}
\bibliographystyle{abbrv}

\end{document}